\documentclass[letterpaper, 10 pt, conference]{ieeeconf}  

\usepackage{amsmath}
\usepackage{amssymb}
\usepackage{mathrsfs}
\usepackage{microtype}
\usepackage{subcaption}

\usepackage{pgfplots}
 \pgfplotsset{compat=newest}
  \usetikzlibrary{plotmarks}
  \usepackage{grffile}
\usepackage{graphics,xcolor}
\definecolor{P285U}{cmyk}{0.89,0.43,0.0,0.0}
\definecolor{P285U_font}{cmyk}{0.89,0.43,0.0,0.4}
\definecolor{lgray}{cmyk}{0,0,0,0.2}
\definecolor{myblue}{cmyk}{100,75,0,0}
\definecolor{jkuBlue}{RGB}{4,110,152}
\definecolor{jkuBlue}{RGB}{0,120,170}
\definecolor{jkuCyan}{RGB}{100,180,190}
\definecolor{jkuYellow}{RGB}{230,195,35}
\definecolor{jkuGrey}{RGB}{125,130,140}
\definecolor{jkuDarkGrey}{RGB}{51,51,51}
\definecolor{jkuLightGreen}{RGB}{195,215,75}
\definecolor{jkuGreen}{RGB}{115,180,85}
\definecolor{jkuPurple}{RGB}{145,75,130}
\definecolor{jkuRed}{RGB}{205,90,80}

\usepackage{mdframed}

\IEEEoverridecommandlockouts                              

\title{\LARGE \bf
Stability Analysis of the Observer Error of an In-Domain Actuated Vibrating String*
}

\author{Tobias Malzer$^{1}$, Hubert Rams$^{2}$, Bernd Kolar$^{1}$ and Markus Sch\"{o}berl$^{1}$
\thanks{*This work has been supported by the Austrian Science Fund (FWF) under
grant number P 29964-N32.}
\thanks{$^{1}$Tobias Malzer, Bernd Kolar and Markus Sch\"{o}berl are with the Institute of Automatic
Control and Control Systems Technology, Johannes Kepler University Linz,
Altenbergerstrasse 66, 4040 Linz, Austria.
        {\tt\small \{tobias.malzer\_1, bernd.kolar, markus.schoeberl\}@jku.at}}%
\thanks{$^{2}$Hubert Rams is with B\&R Industrial Automation GmbH, B\&R Stra\ss{}e 1, 5142 Eggelsberg, Austria.
       {\tt\small hubert.rams@br-automation.com}}%
}

\begin{document}

\maketitle
\thispagestyle{empty}
\pagestyle{empty}

\begin{abstract}

In this paper, the behaviour of the observer error of an in-domain actuated vibrating string, where the observer system has been designed based on energy considerations exploiting a port-Hamiltonian system representation for infinite-dimensional systems, is analysed. Thus, the observer-error dynamics are reformulated as an abstract Cauchy problem, which enables to draw conclusions regarding the well-posedness of the observer-error system. Furthermore, we show that the observer error is asymptotically stable by applying LaSalle's invariance principle for infinite-dimensional systems.

\end{abstract}

\section{Introduction}

A very popular research discipline in control theory is the extension
of control methodologies originally developed for systems that are
described by ordinary differential equations (ODEs) to systems governed
by partial differential equations (PDEs); however, with regard to
stability investigations this extension is accompanied by a significant
rise of complexity, see e.g. \cite{Luo1998} for a comprehensive framework
for the stability analysis of infinite-dimensional systems. Therefore,
a lot of research effort has been invested in this topic, where for
example the stability analysis of mechanical systems with certain
boundary conditions has been addressed. For instance, in \cite{Miletic2015}
the stability of an Euler-Bernoulli beam subjected to nonlinear damping
and a nonlinear spring at the tip is analysed, whereas \cite{Stuerzer2016}
is concerned with the stability behaviour of a gantry crane with heavy
chain and payload. Furthermore, the proof of stability of a Lyapunov-based
control law as well as a Lyapunov-based observer design for an in-domain
actuated Euler-Bernoulli beam has been presented in \cite{Henikl2012}.

A well-known methodology, that has also been extended to the infinite-dimensional
scenario, is the combination of the port-Hamiltonian (pH) system representation
with energy-based control. In this regard, in particular a pH-system
representation based on an underlying jet-bundle structure, see e.g.
\cite{Ennsbrunner2005,Schoeberl2012,Schoeberl2014a}, as well as a
formulation exploiting Stokes-Dirac structures, see e.g. \cite{Schaft2002,Gorrec2005},
have turned out to be especially suitable. For a detailed comparison
of these approaches, where the main difference is the choice of the
variables, the interested reader is referred to \cite{Schoeberl2013b}
or \cite{Malzer2020}. In fact, with respect to boundary-control systems,
a lot of literature is available, see e.g. \cite{Schoeberl2011,Rams2017a}
and \cite{Macchelli2004,Macchelli2017}, where boundary controllers
based on the well-known energy-Casimir method are designed within
the jet-bundle and the Stokes-Dirac approach, respectively. Moreover,
recently the pH-system description has also been exploited with regard
to the observer design, see e.g. \cite{Toledo2019}, where a pH-based
observer-design procedure for boundary-control systems has been developed
within the Stokes-Dirac scenario. In light of the observer design,
of course stability investigations play an important role, since it
must be ensured that the observer error tends to zero.

In \cite{Malzer2020}, a control-design procedure based on the energy-Casimir
method together with an observer design exploiting the pH-system representation
has been presented within the jet-bundle framework as well as within
the Stokes-Dirac scenario for infinite-dimensional systems with in-domain
actuation. Furthermore, the design procedures have been demonstrated
and compared by means of an in-domain actuated vibrating string; however,
the investigation regarding the asymptotic stability of the observer
error -- which is of course essential -- has only been sketched.
Therefore, the aim of this paper is to carry out the stability investigation
of the observer error of this system in detail. To this end, first
of all, in Section \ref{sec:Observer_Design} we summarise the observer
design that exploits the pH-system representation based on a jet-bundle
structure, while in Section \ref{sec:Vibrating_String} the observer
design is explicitly demonstrated for an in-domain actuated vibrating
string. Thus, the main contribution of this paper is to verify the
asymptotic stability of the observer error, where i) it is necessary
to investigate the well-posedness, see Subsection \ref{subsec:Wellposedness},
and ii) to apply LaSalle's invariance principle for infinite-dimensional
systems, see Subsection \ref{subsec:LaSalle}.

\section{Observer Design based on a Port-Hamiltonian Framework\label{sec:Observer_Design}}

With respect to the observer design, see \cite[Sec. V]{Malzer2020},
we intend to exploit a pH-system description for infinite-dimensional
systems with $1$-dimensional spatial domain, which is equipped with
the spatial coordinate $z\in[0,L]$. The system representation is
based on an underlying jet-bundle structure, and therefore, first
of all we introduce the bundle $\pi:\mathcal{E}\rightarrow\mathcal{B}$,
where the total manifold $\mathcal{E}$ is equipped with the coordinates
$(z,x^{\alpha})$, with $x^{\alpha}$, $\alpha=1,\ldots,n$, denoting
the dependent variables, while the base manifold $\mathcal{B}$ possesses
the independent (spatial) coordinate $(z)$ solely. Next, we consider
the so-called vertical tangent bundle $\nu_{\mathcal{E}}:\mathcal{V}(\mathcal{E})\rightarrow\mathcal{E}$,
equipped with the coordinates $(z,x^{\alpha},\dot{x}^{\alpha})$,
which is a subbundle of the tangent bundle $\tau_{\mathcal{E}}:\mathcal{T}(\mathcal{E})\rightarrow\mathcal{E}$,
possessing the coordinates $(z,x^{\alpha},\dot{z},\dot{x}^{\alpha})$
together with the fibre bases $\partial_{z}=\partial/\partial z$
and $\partial_{\alpha}=\partial/\partial x^{\alpha}$. Thus, a vertical
vector field $v=\mathcal{E}\rightarrow\mathcal{V}(\mathcal{E})$,
in local coordinates given as $v=v^{\alpha}\partial_{\alpha}$ with
$v^{\alpha}\in C^{\infty}(\mathcal{E})$, i.e. $v^{\alpha}$ is a
smooth function on $\mathcal{E}$, is defined as a section. A further
important differential-geometric object is the so-called co-tangent
bundle $\tau_{\mathcal{E}}^{*}:\mathcal{T}^{*}(\mathcal{E})\rightarrow\mathcal{E}$,
possessing the coordinates $(z,x^{\alpha},\dot{z},\dot{x}_{\alpha})$
together with the fibre bases $\mathrm{d}z$ and $\mathrm{d}x^{\alpha}$,
which allows to introduce a one-form $w:\mathcal{E}\rightarrow\mathcal{T}^{*}(\mathcal{E})$
as a section that can locally be given as $w=\breve{w}\mathrm{d}z+w_{\alpha}\mathrm{d}x^{\alpha}$
with $\breve{w},w_{\alpha}\in C^{\infty}(\mathcal{E})$. With respect
to the pH-system representation, we are interested in densities $\mathfrak{H}=\mathcal{H}\mathrm{d}z$
with $\mathcal{H}\in C^{\infty}(\mathcal{J}^{1}(\mathcal{E}))$, where
these densities can be formed by sections of certain pullback bundles,
whose use is omitted here for ease of presentation. That is, $\mathcal{H}$
is a smooth function on the first jet manifold $\mathcal{J}^{1}(\mathcal{E})$,
which is equipped with the coordinates $(z,x^{\alpha},x_{z}^{\alpha})$,
where the $1$st-order jet variable $x_{z}^{\alpha}$ corresponds
to the derivative of $x^{\alpha}$ with respect to $z$. Moreover,
the first prolongation of a vertical vector field reads as $j^{1}(v)=v^{\alpha}\partial_{\alpha}+d_{z}(v^{\alpha})\partial_{\alpha}^{z}$,
with $\partial_{\alpha}^{z}=\partial/\partial x_{z}^{\alpha}$, where
we exploit the total derivative $d_{z}=\partial_{z}+x_{z}^{\alpha}\partial_{\alpha}+x_{zz}^{\alpha}\partial_{\alpha}^{z}+\ldots$.

Having discussed this essential preliminaries, we are able to introduce
the pH-system representation including inputs and outputs on the spatial
domain as\begin{subequations}\label{eq:pH_sys_jetbundle}
\begin{align}
\dot{x} & =(\mathcal{J}-\mathcal{R})(\delta\mathfrak{H})+u\rfloor\mathcal{G}\,,\label{eq:pH_sys_dynamics}\\
y & =\mathcal{G}^{*}\rfloor\delta\mathfrak{H}\,,
\end{align}
\end{subequations}see e.g. \cite{Ennsbrunner2005,Schoeberl2008a,Schoeberl2014},
where $\rfloor$ denotes the so-called Hook operator allowing for
the natural contraction between tensor fields. In (\ref{eq:pH_sys_jetbundle}),
the variational derivative $\delta\mathfrak{H}=\delta_{\alpha}\mathcal{H}\mathrm{d}x^{\alpha}\wedge\mathrm{d}z$,
with $\wedge$ denoting the exterior (wedge) product, locally reads
as $\delta_{\alpha}\mathcal{H}=\partial_{\alpha}\mathcal{H}-d_{z}(\partial_{\alpha}^{z}\mathcal{H})$.
Furthermore, the linear operators $\mathcal{J},\mathcal{R}:\mathcal{T}^{*}(\mathcal{E})\wedge\mathcal{T}^{*}(\mathcal{B})\rightarrow\mathcal{V}(\mathcal{E})$
describe the internal power flow and the dissipation effects of the
system, respectively. The coefficients $\mathcal{J}^{\alpha\beta}$
of the interconnection tensor $\mathcal{J}$ meet $\mathcal{J}^{\alpha\beta}=-\mathcal{J}^{\beta\alpha}\in C^{\infty}(\mathcal{J}^{2}(\mathcal{E}))$,
while we have $\mathcal{R}^{\alpha\beta}=\mathcal{R}^{\beta\alpha}\in C^{\infty}(\mathcal{J}^{2}(\mathcal{E}))$
and $[\mathcal{R}^{\alpha\beta}]\geq0$ for the coefficient matrix
of the symmetric and positive semi-definite dissipation mapping $\mathcal{R}$.
With respect to the dual input and output bundles $\rho:\mathcal{U}\rightarrow\mathcal{J}^{2}(\mathcal{E})$
and $\varrho:\mathcal{Y}\rightarrow\mathcal{J}^{2}(\mathcal{E})$,
we have the input map and its adjoint output map $\mathcal{G}:\mathcal{U}\rightarrow\mathcal{V}(\mathcal{E})$
and $\mathcal{G}^{*}:\mathcal{T}^{*}(\mathcal{E})\wedge\mathcal{T}^{*}(\mathcal{B})\rightarrow\mathcal{Y}$,
respectively, and thus, the relation $(u\rfloor\mathcal{G})\rfloor\delta\mathfrak{H}=u\rfloor(\mathcal{G}^{*}\rfloor\delta\mathfrak{H})=u\rfloor y$
holds, see \cite[Sec. 4]{Ennsbrunner2005} or \cite[Sec. 3]{Schoeberl2008a}.
To be able to determine the formal change of the Hamiltonian functional
$\mathscr{H}=\int_{0}^{L}\mathcal{H}\mathrm{d}z$ along solutions
of (\ref{eq:pH_sys_dynamics}), we make use of the Lie-derivative
$\mathrm{L}_{j^{1}(v)}$, where we set $v=\dot{x}$ with (\ref{eq:pH_sys_dynamics}),
see \cite[Sec. IV-A]{Schoeberl2011}, and thus, we obtain
\begin{equation}
\dot{\mathscr{H}}=-\int_{0}^{L}\mathcal{R}(\delta\mathfrak{H})\rfloor\delta\mathfrak{H}+\int_{0}^{L}u\rfloor y+\left.(\dot{x}\rfloor\delta^{\partial}\mathfrak{H})\right|_{0}^{L}\,\label{eq:H_p}
\end{equation}
by means of integration by parts and Stoke's theorem. If $\mathscr{H}$
corresponds to the total energy of the system, then (\ref{eq:H_p})
states a power-balance relation, where the first expression describes
the energy that is dissipated for example due to damping effects.
Moreover, the expression $\int_{0}^{L}u\rfloor y$ denotes a collocation
term distributed over (a part of) the spatial domain. The last term
corresponds to collocation restricted to the boundary, which is indicated
by $(\cdot)|_{0}^{L}$, where the boundary operator locally reads
as $\delta_{\alpha}^{\partial}\mathcal{H}=\partial_{\alpha}^{z}\mathcal{H}$.
Note that here we consider systems with trivial boundary conditions,
implying that the boundary ports $(\dot{x}^{\alpha}\delta_{\alpha}^{\partial}\mathcal{H})|_{0}^{L}$
vanish.

Next, the intention is to exploit the pH-formulation with respect
to the observer design. In particular, the copy of the plant (\ref{eq:pH_sys_dynamics})
is extended by an error-injection term, and thus, by means of the
observer-energy density $\hat{\mathcal{H}}$, the observer system
is locally given by\begin{subequations}\label{eq:observer_JB}
\begin{align}
\dot{\hat{x}}^{\hat{\alpha}} & =(\mathcal{J}^{\hat{\alpha}\hat{\beta}}-\mathcal{R}^{\hat{\alpha}\hat{\beta}})\delta_{\hat{\beta}}\hat{\mathcal{H}}+\mathcal{G}_{\xi}^{\hat{\alpha}}u^{\xi}+\mathcal{K}_{\eta}^{\hat{\alpha}}u_{o}^{\eta}\,,\label{eq:observer_JB_dynamics}\\
\hat{y}_{\xi} & =\mathcal{G}_{\xi}^{\hat{\alpha}}\delta_{\hat{\alpha}}\hat{\mathcal{H}}\,,\label{eq:observer_JB_output_densities}
\end{align}
\end{subequations}with $\hat{\alpha},\hat{\beta}=1,\ldots,n$ and
$\xi,\eta=1,\ldots,m$, where we use Einstein's convention on sums.
In (\ref{eq:observer_JB_dynamics}), we have the additional input
$u_{o}^{\eta}=\delta^{\eta\xi}(\bar{y}_{\xi}-\hat{\bar{y}}_{\xi})$
-- with the Kronecker-Delta symbol meeting $\delta^{\xi\eta}=1$
for $\xi=\eta$ and $\delta^{\xi\eta}=0$ for $\xi\neq\eta$ --,
where $\bar{y}_{\xi}$ corresponds to the integrated output density
of the plant according to $\bar{y}_{\xi}=\int_{0}^{L}y_{\xi}\mathrm{d}z$,
which is assumed to be available as measurement quantity, while $\hat{\bar{y}}_{\xi}$
represents the copy of the integrated plant-output according to $\hat{\bar{y}}_{\xi}=\int_{0}^{L}\hat{y}_{\xi}\mathrm{d}z$
with (\ref{eq:observer_JB_output_densities}). The aim is to design
the observer gain $\mathcal{K}_{\eta}^{\hat{\alpha}}$ such that the
observer error $\tilde{x}=x-\hat{x}$ tends to $0$, where it is beneficial
to reformulate the observer-error dynamics $\dot{\tilde{x}}=\dot{x}-\dot{\hat{x}}$
as pH-system according to\begin{subequations}\label{eq:observer_error_system_JB}
\begin{align}
\dot{\tilde{x}}^{\tilde{\alpha}} & =(\mathcal{J}^{\tilde{\alpha}\tilde{\beta}}-\mathcal{R}^{\tilde{\alpha}\tilde{\beta}})\delta_{\tilde{\beta}}\tilde{\mathcal{H}}-\mathcal{K}_{\xi}^{\tilde{\alpha}}u_{o}^{\xi}\,,\\
\tilde{y}_{\xi} & =-\mathcal{K}_{\xi}^{\tilde{\alpha}}\delta_{\tilde{\alpha}}\tilde{\mathcal{H}}\,.\label{eq:observer_error_output_JB}
\end{align}
\end{subequations}with (\ref{eq:observer_error_output_JB}) denoting
the collocated output density. If we investigate the formal change
of the error-Hamiltonian $\tilde{\mathscr{H}}=\int_{0}^{L}\tilde{\mathcal{H}}\mathrm{d}z$,
which follows to
\begin{multline*}
\dot{\tilde{\mathscr{H}}}=-\int_{0}^{L}\delta_{\tilde{\alpha}}(\tilde{\mathcal{H}})\mathcal{R}^{\tilde{\alpha}\tilde{\beta}}\delta_{\tilde{\beta}}(\tilde{\mathcal{H}})\mathrm{d}z+\ldots\\
-\int_{0}^{L}\delta_{\tilde{\alpha}}(\tilde{\mathcal{H}})\mathcal{K}_{\xi}^{\tilde{\alpha}}\delta^{\xi\eta}(\bar{y}_{\eta}-\hat{\bar{y}}_{\eta})\mathrm{d}z\,,
\end{multline*}
we find that by means of a proper choice for the components $\mathcal{K}_{\xi}^{\hat{\alpha}}$
we are able to render $\dot{\tilde{\mathscr{H}}}\leq0$ . Hence, the
total energy of the observer error $\tilde{\mathscr{H}}$ is an appropriate
candidate for a Lyapunov functional and therefore serves as basis
with respect to the stability analysis. Next, the observer-design
procedure is demonstrated by an example.

\section{Observer Design for an In-Domain Actuated Vibrating String\label{sec:Vibrating_String}}

In this chapter, we design an infinite-dimensional observer for an
in-domain actuated vibrating string by exploiting energy considerations.
The governing equation of motion of the system under consideration
reads as\begin{subequations}
\begin{equation}
\rho\frac{\partial^{2}w}{\partial t^{2}}=T\frac{\partial^{2}w}{\partial z^{2}}+f(z,t)\,,\label{eq:vib_String_eom}
\end{equation}
where $w$ describes the vertical deflection of the string, $\rho$
the mass density and $T$ Young's modulus. Regarding the boundary
conditions, we have that the string is clamped at $z=0$ and free
at $z=L$, i.e.
\begin{align}
w(0,t) & =0\,,\quad T\frac{\partial w}{\partial z}(L,t)=0\,.\label{eq:BC_VS}
\end{align}
\end{subequations}In (\ref{eq:vib_String_eom}), the distributed
force $f(z,t)=g(z)u(t)$ is generated by an actuator behaving like
a piezoelectric patch, where the applied voltage $u(t)$ serves as
manipulated variable. The spatially dependent function $g(z)=h(z-L_{p_{1}})-h(z-L_{p_{2}})$,
where $h(\cdot)$ denotes the Heaviside function, describes the placement
of the actuator between $z=L_{p_{1}}$ and $z=L_{p_{2}}$. In fact,
the force-distribution on the domain $L_{p_{1}}\leq z\leq L_{p_{2}}$
is supposed to be constant and is scaled by $u(t)$.

First, the intention is to find a pH-system representation that can
be exploited for the observer design. To this end, we introduce the
underlying bundle structure based on $\pi:(z,w,p)\rightarrow(z)$
together with the generalised momenta $p=\rho\dot{w}$, and thus,
(\ref{eq:vib_String_eom}) can be rewritten as
\begin{equation}
\dot{p}=Tw_{zz}+g(z)u\,.\label{eq:VS_JB}
\end{equation}
If we use the Hamiltonian density $\mathcal{H}=\frac{1}{2\rho}p^{2}+\frac{1}{2}T(w_{z})^{2}\in\mathcal{J}^{1}(\mathcal{E})$,
we obtain the appropriate pH-system formulation\begin{subequations}\label{eq:pH_formulation_VS}
\begin{align}
\left[\begin{array}{c}
\dot{w}\\
\dot{p}
\end{array}\right] & =\left[\begin{array}{cc}
0 & 1\\
-1 & 0
\end{array}\right]\left[\begin{array}{c}
\delta_{w}\mathcal{H}\\
\delta_{p}\mathcal{H}
\end{array}\right]+\left[\begin{array}{c}
0\\
g(z)
\end{array}\right]u\,,\\
y & =\left[\begin{array}{cc}
0 & g(z)\end{array}\right]\left[\begin{array}{c}
\delta_{w}\mathcal{H}\\
\delta_{p}\mathcal{H}
\end{array}\right]=g(z)\frac{p}{\rho}\,.\label{eq:pH_form_VS_output_density}
\end{align}
\end{subequations}By taking the boundary conditions (\ref{eq:BC_VS})
into account, one finds that the formal change of the Hamiltonian
functional $\mathscr{H}$ follows to $\dot{\mathscr{H}}=\int_{0}^{L}g(z)\frac{p}{\rho}u\mathrm{d}z$,
i.e. we have a distributed port that can be used for control purposes.
In fact, for the system under consideration, in \cite{Malzer2020}
a dynamic controller based on the energy-Casimir method has been designed.
However, with regard to this control methodology, it should be mentioned
that it yields unsatisfactory results for uncertain initial conditions,
see e.g. \cite{Rams2017b} where this problem is briefly discussed
for a boundary-control system. Therefore, in the following we intend
to design an infinite-dimensional observer in order to overcome this
obstacle.

Concerning the observer design, it is assumed that the spatial integration
of the distributed output density (\ref{eq:pH_form_VS_output_density})
according to $\bar{y}=\int_{0}^{L}g(z)\frac{p}{\rho}\mathrm{d}z$,
which can be interpreted as the current through the actuator, is available
as measurement quantity. Thus, if we use the observer density $\hat{\mathcal{H}}=\frac{1}{2\rho}\hat{p}^{2}+\frac{1}{2}T(\hat{w}_{z})^{2}$
and the copy of the plant output according to $\hat{\bar{y}}=\int_{0}^{L}g(z)\frac{\hat{p}}{\rho}\mathrm{d}z$,
we are able to introduce an observer for the in-domain actuated vibrating
string in the form
\begin{align*}
\left[\begin{array}{c}
\dot{\hat{w}}\\
\dot{\hat{p}}
\end{array}\right]\! & =\!\left[\!\begin{array}{cc}
0\! & \!1\\
-1\! & \!0
\end{array}\!\right]\!\left[\begin{array}{c}
\delta_{\hat{w}}\hat{\mathcal{H}}\\
\delta_{\hat{p}}\hat{\mathcal{H}}
\end{array}\right]\!+\!\left[\begin{array}{c}
0\\
g(z)
\end{array}\right]\!u\!+\!\left[\begin{array}{c}
k_{1}\\
k_{2}
\end{array}\right]\!(\bar{y}-\hat{\bar{y}}),
\end{align*}
where the governing equations are restricted to the boundary conditions
$\hat{w}(0)=0$ and $T\hat{w}_{z}(L)=0$. Next, by means of the error
coordinates $\tilde{w}=w-\hat{w}$, $\tilde{p}=p-\hat{p}$, the observer-error
dynamics can be deduced to\begin{subequations}\label{eq:observer_error_dynamics_VS}
\begin{align}
\dot{\tilde{w}} & =\dot{w}-\dot{\hat{w}}=\frac{1}{\rho}\tilde{p}-k_{1}(\bar{y}-\hat{\bar{y}})\,,\label{eq:observer_error_dyn_w}\\
\dot{\tilde{p}} & =\dot{p}-\dot{\hat{p}}=T\tilde{w}_{zz}-k_{2}(\bar{y}-\hat{\bar{y}})\,,\label{eq:observer_error_dyn_p}
\end{align}
where the boundary conditions
\begin{align}
\tilde{w}(0) & =0\,,\quad T\tilde{w}_{z}(L)=0\label{eq:BC_observer_error_left}
\end{align}
\end{subequations}hold. With respect to the determination of $k_{1}$
and $k_{2}$, it is beneficial to reformulate (\ref{eq:observer_error_dyn_w})
and (\ref{eq:observer_error_dyn_p}) as the pH-system\begin{subequations}\label{eq:observer_error_VS_JB-1}
\begin{align}
\left[\begin{array}{c}
\dot{\tilde{w}}\\
\dot{\tilde{p}}
\end{array}\right]\! & =\!\left[\begin{array}{cc}
0 & 1\\
-1 & 0
\end{array}\right]\!\left[\begin{array}{c}
\delta_{\tilde{w}}\tilde{\mathcal{H}}\\
\delta_{\tilde{p}}\tilde{\mathcal{H}}
\end{array}\right]\!-\!\left[\begin{array}{c}
k_{1}\\
k_{2}
\end{array}\right]\!(\bar{y}-\hat{\bar{y}}),\\
\tilde{y} & \!=\!-\!\left[\begin{array}{cc}
k_{1} & k_{2}\end{array}\right]\!\left[\begin{array}{c}
\delta_{\tilde{w}}\tilde{\mathcal{H}}\\
\delta_{\tilde{p}}\tilde{\mathcal{H}}
\end{array}\right]\!=\!k_{1}T\tilde{w}_{zz}\!-\!k_{2}\frac{\tilde{p}}{\rho},\label{eq:pH_observer_error_output}
\end{align}
\end{subequations}where the energy density of the observer error
reads as $\tilde{\mathcal{H}}=\frac{1}{2\rho}\tilde{p}^{2}+\frac{1}{2}T(\tilde{w}_{z})^{2}$
and (\ref{eq:pH_observer_error_output}) states the corresponding
output density. If we investigate the formal change of the error-Hamiltonian
functional $\tilde{\mathscr{H}}$, which can be deduced to
\begin{align}
\dot{\tilde{\mathscr{H}}} & =\int_{0}^{L}(T\tilde{w}_{zz}k_{1}(\bar{y}-\hat{\bar{y}})-\frac{\tilde{p}}{\rho}k_{2}(\bar{y}-\hat{\bar{y}}))\mathrm{d}z\,,\label{eq:H_tilde_p}
\end{align}
and take into account that $(\bar{y}-\hat{\bar{y}})=\int_{0}^{L}g(z)\frac{1}{\rho}\tilde{p}\mathrm{d}z$,
we find that the choice $k_{1}=0$ and $k_{2}=kg(z)$ with $k>0$,
yields
\begin{equation}
\dot{\tilde{\mathscr{H}}}(\tilde{w},\tilde{p})=-k(\bar{y}-\hat{\bar{y}})^{2}\leq0\,.\label{eq:H_tilde_p_fin}
\end{equation}
However, the fact that $\tilde{\mathscr{H}}>0$ and $\dot{\tilde{\mathscr{H}}}\leq0$
hold is not sufficient for the convergence of the observer, and therefore,
in the following, detailed stability investigations are carried out
to verify that the observer error is asymptotically stable.

\section{Observer Convergence\label{sec:Stability_Analysis}}

In this section, based on functional analysis the convergence of the
observer error is proven in two steps. First, we address the well-posedness
of the observer-error system making heavy use of the well-known Lumer-Phillips
theorem, see e.g. \cite{Liu1999}. Afterwards, LaSalle's invariance
principle for infinite-dimensional systems is applied to show the
asymptotic stability of the observer error, where beforehand it is
necessary to verify the precompactness of the solution trajectories.

\subsection{Well-posedness of the Observer-Error System\label{subsec:Wellposedness}}

Now, a careful investigation of the well-posedness of the observer-error
system is carried out. To this end, we reformulate (\ref{eq:observer_error_dynamics_VS})
as an abstract Cauchy problem and show that the operator under consideration
generates a $C_{0}$-semigroup of contractions.

First, we define the state vector $\chi=\left[\chi^{1},\chi^{2}\right]^{T}=\left[\tilde{w},\tilde{p}\right]^{T}$
together with the state space $\mathcal{X}=H_{C}^{1}(0,L)\times L^{2}(0,L)$,
where $H_{C}^{1}(0,L)=\{\chi^{1}\in H^{1}(0,L)|\chi^{1}(0)=0\}$,
with $H^{l}(0,L)$ denoting a Sobolev space of functions whose derivatives
up to order $l$ are square integrable, see \cite{Adams2003} for
a detailed introduction of Sobolev spaces. Thus, the state space $\mathcal{X}$
is equipped with the standard norm
\begin{equation}
\left\Vert \chi\right\Vert _{n}^{2}=\left\langle \tilde{w},\tilde{w}\right\rangle _{L^{2}}+\left\langle \tilde{w}_{z},\tilde{w}_{z}\right\rangle _{L^{2}}+\left\langle \tilde{p},\tilde{p}\right\rangle _{L^{2}}\,.\label{eq:standard_norm}
\end{equation}
Next, to be able to rewrite the observer-error dynamics as an abstract
Cauchy problem of the form $\dot{\chi}(t)=\mathcal{A}\chi(t)$ with
$\chi(0)=\chi_{0}$, we introduce the linear operator $\mathcal{A}:\mathcal{D}(\mathcal{A})\subset\mathcal{X}\rightarrow\mathcal{X}$
according to
\[
\mathcal{A}:\left[\begin{array}{c}
\tilde{w}\\
\tilde{p}
\end{array}\right]\rightarrow\left[\begin{array}{c}
\frac{1}{\rho}\tilde{p}\\
T\tilde{w}_{zz}-kg(z)\int_{0}^{L}g(z)\frac{1}{\rho}\tilde{p}\mathrm{d}z
\end{array}\right]\,,
\]
where the (dense) domain of $\mathcal{A}$ is defined as
\begin{multline}
\mathcal{D}(\mathcal{A}):=\{\chi\in\mathcal{X}|\tilde{w}\in(H^{2}(0,L)\cap H_{C}^{1}(0,L)),\\
\tilde{p}\in H_{C}^{1}(0,L),T\tilde{w}_{z}(L)=0\}\,.\label{eq:domain_A}
\end{multline}
Thus, the intention is to investigate the operator $\mathcal{A}$
regarding some properties such that a variant of the well-known Lumer-Phillips
theorem \cite[Thm. 1.2.4]{Liu1999} can be applied. With respect to
this forthcoming investigations, it is beneficial to introduce
\begin{equation}
\left\Vert \chi\right\Vert _{\mathcal{\mathcal{X}}}^{2}=T\left\langle \tilde{w}_{z},\tilde{w}_{z}\right\rangle _{L^{2}}+\frac{1}{\rho}\left\langle \tilde{p},\tilde{p}\right\rangle _{L^{2}}\,,\label{eq:energy_norm}
\end{equation}
which is called energy norm due to the equivalence $\tilde{\mathscr{H}}=\frac{1}{2}\left\Vert \chi\right\Vert _{\mathcal{\mathcal{X}}}^{2}$.
Because $\tilde{w}(0)=0$ and further $\tilde{w}(z)=\int_{0}^{z}\tilde{w}_{z}\mathrm{d}y_{1}$
holds, we find constants $c_{1},c_{2}$, which have to meet $0<c_{1}\leq\mathrm{min}\left(\frac{T}{L+1},\frac{1}{\rho}\right)$
and $c_{2}\geq\mathrm{max}\left(T,\frac{1}{\rho}\right)>0$, such
that $c_{1}\left\Vert \chi\right\Vert _{n}^{2}\leq\left\Vert \chi\right\Vert _{\mathcal{X}}^{2}\leq c_{2}\left\Vert \chi\right\Vert _{n}^{2}$
is fulfilled, and hence, the energy norm (\ref{eq:energy_norm}) is
equivalent to the standard norm (\ref{eq:standard_norm}). Similar
to the proof of Lemma 2.2 in \cite{Stuerzer2016}, where they exploit
the dense inclusion $H^{2}(0,L)\subset H^{1}(0,L)$ and modify the
boundary values of $\tilde{w}$ and its derivatives in a proper manner,
it can be shown that the domain $\mathcal{D}(\mathcal{A})$ given
in (\ref{eq:domain_A}) is dense in $\mathcal{X}$. Thus, according
to \cite[Def. 1.1.1]{Liu1999}, -- since we have the equivalence
$\tilde{\mathscr{H}}=\frac{1}{2}\left\Vert \chi\right\Vert _{\mathcal{\mathcal{X}}}^{2}$
-- the relation (\ref{eq:H_tilde_p_fin}) implies that $\mathcal{A}$
is dissipative.

In the following, we show that the inverse operator $\mathcal{A}^{-1}$
exists and is bounded, i.e. for every $\bar{\chi}=\left[f,h\right]^{T}\in\mathcal{X}$
and $\chi=\left[\tilde{w},\tilde{p}\right]^{T}\in\mathcal{D}(\mathcal{A})$,
we can uniquely solve
\begin{equation}
\mathcal{A}\!\left[\!\begin{array}{c}
\tilde{w}\\
\tilde{p}
\end{array}\!\right]\!=\!\left[\begin{array}{c}
\frac{1}{\rho}\tilde{p}\\
T\tilde{w}_{zz}-kg(z)\int_{0}^{L}g(z)\frac{1}{\rho}\tilde{p}\mathrm{d}z
\end{array}\right]\!=\!\left[\begin{array}{c}
f\\
h
\end{array}\right],\label{eq:Calc_A-1}
\end{equation}
and prove that $\mathcal{A}^{-1}$ maps bounded sets in $\mathcal{X}$
into bounded sets in $\mathcal{K}:=(H^{2}(0,L)\cap H_{C}^{1}(0,L))\times H_{C}^{1}(0,L)$.
From the $1$st line of (\ref{eq:Calc_A-1}) it follows that $\tilde{p}=\rho f\in H_{C}^{1}(0,L)$.
Moreover, an integration of the $2$nd line of (\ref{eq:Calc_A-1})
yields
\begin{multline}
\tilde{w}_{z}(z)=-\frac{1}{T}(\int_{z}^{L}h(y_{2})\mathrm{d}y_{2}+\ldots\\
+\int_{z}^{L}kg(y_{2})\int_{0}^{L}g(y_{1})f(y_{1})\mathrm{d}y_{1}\mathrm{d}y_{2})\label{eq:w_z_calc}
\end{multline}
as $\tilde{w}_{z}(L)=0$ holds. If we further integrate (\ref{eq:w_z_calc}),
we obtain
\begin{multline}
\tilde{w}(z)=-\frac{1}{T}(\int_{0}^{z}\int_{y_{3}}^{L}h(y_{2})\mathrm{d}y_{2}\mathrm{d}y_{3}+\ldots\\
+\int_{0}^{z}\int_{y_{3}}^{L}kg(y_{2})\int_{0}^{L}g(y_{1})f(y_{1})\mathrm{d}y_{1}\mathrm{d}y_{2}\mathrm{d}y_{3})\label{eq:w_calc}
\end{multline}
as $\tilde{w}(0)=0$, and thus, $\tilde{w}(z)$ is uniquely defined
by $\bar{\chi}$. Since we have shown that the inverse operator $\mathcal{A}^{-1}$
exists, it remains to investigate the boundedness. To this end, it
is verified that the norm of $\chi=\mathcal{A}^{-1}\bar{\chi}$ in
$\mathcal{K}$ is bounded by $\left\Vert \bar{\chi}\right\Vert _{\mathcal{X}}$.
First, we state an inequality that is often used in the sequel; in
fact, for a -- basically arbitrary -- function $f$, by means of
the Cauchy-Schwarz inequality we find the important relation
\begin{equation}
(\int_{0}^{L}f\mathrm{d}z)^{2}\leq C\int_{0}^{L}\left|f\right|^{2}\mathrm{d}z\,,\label{eq:inequalitiy_square}
\end{equation}
where it should be mentioned that here and in the following $C$
denotes positive, not necessarily equal constants. Next, we investigate
the norm $\left\Vert \tilde{w}_{z}\right\Vert _{L^{2}}$. Therefore,
we substitute (\ref{eq:w_z_calc}) in $\left\Vert \tilde{w}_{z}\right\Vert _{L^{2}}=(\int_{0}^{L}\left|\tilde{w}_{z}\right|^{2}\mathrm{d}z)^{1/2}$
and apply the Triangle inequality, which yields
\begin{multline*}
\left\Vert \tilde{w}_{z}\right\Vert _{L^{2}}\leq(\int_{0}^{L}\frac{1}{T^{2}}(\int_{z}^{L}h(y_{2})\mathrm{d}y_{2})^{2}\mathrm{d}z)^{\frac{1}{2}}+\\
(\int_{0}^{L}\frac{1}{T^{2}}(\int_{z}^{L}kg(y_{2})\int_{0}^{L}g(y_{1})f(y_{1})\mathrm{d}y_{1}\mathrm{d}y_{2})^{2}\mathrm{d}z)^{\frac{1}{2}}.
\end{multline*}
 Thus, by means of (\ref{eq:inequalitiy_square}) and due to the fact
that $\int_{z}^{L}h^{2}\mathrm{d}y_{2}\leq\int_{0}^{L}h^{2}\mathrm{d}z=\left\Vert h\right\Vert _{L_{2}}^{2}$
holds, we obtain
\begin{multline}
\left\Vert \tilde{w}_{z}\right\Vert _{L^{2}}\leq C\left\Vert h\right\Vert _{L_{2}}\frac{1}{T}L^{\frac{1}{2}}+\\
C(\int_{0}^{L}\frac{1}{T^{2}}\int_{z}^{L}k^{2}g^{2}(y_{2})(\int_{0}^{L}g(y_{1})f(y_{1})\mathrm{d}y_{1})^{2}\mathrm{d}y_{2}\mathrm{d}z)^{\frac{1}{2}}.\label{eq:ineq_w_z_L_2}
\end{multline}
Next, we apply the Cauchy-Schwarz inequality to the second term of
the right-hand side in (\ref{eq:ineq_w_z_L_2}), which enables us
to find the estimate $\left\Vert \tilde{w}_{z}\right\Vert _{L^{2}}\leq C(\left\Vert f\right\Vert _{H^{1}}+\left\Vert h\right\Vert _{L^{2}})$.
Similarly, by means of the $2$nd line of (\ref{eq:Calc_A-1}) we
are able to deduce $\left\Vert \tilde{w}_{zz}\right\Vert _{L^{2}}\leq C(\left\Vert f\right\Vert _{H^{1}}+\left\Vert h\right\Vert _{L^{2}})$.
Moreover, if we substitute (\ref{eq:w_calc}) in $\left\Vert \tilde{w}\right\Vert _{L^{2}}=(\int_{0}^{L}\left|\tilde{w}\right|^{2}\mathrm{d}z)^{1/2}$,
we find $\left\Vert \tilde{w}\right\Vert _{L^{2}}\leq C(\left\Vert f\right\Vert _{H^{1}}+\left\Vert h\right\Vert _{L^{2}})$,
and hence, we have $\left\Vert \tilde{w}\right\Vert _{H^{2}}\leq C(\left\Vert f\right\Vert _{H^{1}}+\left\Vert h\right\Vert _{L^{2}})$.
Since from the first line in (\ref{eq:Calc_A-1}) we immediately get
$\left\Vert \tilde{p}\right\Vert _{H^{1}}=\rho\left\Vert f\right\Vert _{H^{1}}$,
we can state the important estimate
\[
\left\Vert \tilde{w}\right\Vert _{H^{2}}+\left\Vert \tilde{p}\right\Vert _{H^{1}}\leq C(\left\Vert f\right\Vert _{H^{1}}+\left\Vert h\right\Vert _{L^{2}})\,,
\]
which shows that $\mathcal{A}^{-1}$ maps bounded sets in $\mathcal{X}$
into bounded sets in $\mathcal{K}$. 

The boundedness of $\mathcal{A}^{-1}$ implies that $\lambda=0$ cannot
be an eigenvalue of $\mathcal{A}$, and hence, it follows that $0\in\rho(\mathcal{A})$,
the resolvent set of $\mathcal{A}$. Furthermore, since $\mathcal{D}(\mathcal{A})$
is dense in $\mathcal{X}$ and $\mathcal{A}$ is dissipative, all
requirements for the variant of the Lumer-Phillips theorem according
to \cite[Thm. 1.2.4]{Liu1999} are met, and therefore, we are able
to show that $\mathcal{A}$ is the infinitesimal generator of a $C_{0}$-semigroup
of contractions on $\mathcal{X}$. That is, the norm $\left\Vert \chi(t)\right\Vert _{\mathcal{X}}$
remains bounded for $t\rightarrow\infty$; however, with respect to
the observer error it is necessary that it tends to $0$, which is
shown in the following subsection.

\subsection{Asymptotic Stability of the Observer-Error System\label{subsec:LaSalle}}

Now, the objective is to apply LaSalle's invariance principle in order
to prove the asymptotic stability of the observer error, where the
proof follows the intention of \cite[Sec. 3]{Guo2011}. However, the
applicability of LaSalle's invariance principle according to \cite[Thm. 3.64]{Luo1998}
requires the precompactness of the solution trajectories, which is
not ensured in the infinite-dimensional scenario. Since in the previous
section we have shown that $\mathcal{A}^{-1}$ is bounded, by means
of the Sobolev embedding theorem, it follows that $\mathcal{A}^{-1}$
is compact (see proof of Lemma 2.4 in \cite{Stuerzer2016} or \cite[p. 201]{Luo1998}),
which further implies the precompactness of the trajectories, see
\cite[Rem. 4.2]{Miletic2015}. 

In light of LaSalle's invariance principle, we investigate the set
$\mathcal{S}=\{\chi\in\mathcal{X}|\dot{\tilde{\mathscr{H}}}=0\}$,
where $\dot{\tilde{\mathscr{H}}}(\tilde{w},\tilde{p})=-k(\int_{0}^{L}g(z)\frac{1}{\rho}\tilde{p}\mathrm{d}z)^{2}=0$
implies $\int_{0}^{L}g(z)\frac{1}{\rho}\tilde{p}\mathrm{d}z=0$. In
the set $\mathcal{S}$ we have\begin{subequations}\label{eq:eom_S}
\begin{align}
\rho\tilde{w}_{tt} & =T\tilde{w}_{zz}\,,\label{eq:eom_S_pde}\\
\tilde{w}(0,t) & =0\,,\label{eq:eom_S_BC_left}\\
T\tilde{w}_{z}(L,t) & =0\,,\label{eq:eom_S_BC_right}
\end{align}
\end{subequations}which is similar to the problem considered in \cite[Sec. 3]{Guo2011};
however, the restriction describing the set $\mathcal{S}$, which
is constrained to the boundary there, is completely different. To
be able to show that the only possible solution in $\mathcal{S}$
is the trivial one, we need to investigate the general solution of
(\ref{eq:eom_S}). To this end, like in \cite[Sec. 3]{Guo2011}, we
first focus on determining the eigenvalues and eigenfunctions of (\ref{eq:eom_S}),
i.e. we consider
\begin{equation}
\bar{\mathcal{A}}\left[\begin{array}{cc}
\phi & \kappa\end{array}\right]^{T}=\left[\begin{array}{cc}
\frac{\kappa}{\rho} & T\phi_{zz}\end{array}\right]=\lambda\left[\begin{array}{cc}
\phi & \kappa\end{array}\right]^{T}\,.\label{eq:eigen_equation}
\end{equation}
From (\ref{eq:eigen_equation}) we obtain $\kappa=\rho\lambda\phi$
and furthermore\begin{subequations}\label{eq:eigen_problem}
\begin{align}
\phi_{zz} & =\frac{\lambda^{2}}{\vartheta^{2}}\phi\,,\\
\phi(0) & =0\,,\label{eq:eigen_problem_BC_left}\\
\phi_{z}(L) & =0\,,\label{eq:eigen_problem_BC_right}
\end{align}
\end{subequations}where $\vartheta^{2}=\frac{T}{\rho}$. To find
the solution of (\ref{eq:eigen_problem}), we have to investigate
the three cases $\lambda^{2}>0$, $\lambda^{2}=0$ and $\lambda^{2}<0$
in the following. For $\lambda^{2}>0$ and $\lambda^{2}=0$, we have
the ansatz $\phi(z)=Ae^{\frac{\lambda}{\vartheta}z}+Be^{-\frac{\lambda}{\vartheta}z}$
and $\phi(z)=Az+B$, respectively, where by means of the boundary
conditions (\ref{eq:eigen_problem_BC_left}) and (\ref{eq:eigen_problem_BC_right}),
one can easily deduce that for both cases only the trivial solution
$\phi(z)=0$ exists. Thus, we focus on the case $\lambda^{2}<0$,
and consequently, due to the fact that $\lambda$ has an imaginary
character then, as ansatz for the eigenfunctions we have $\phi(z)=A\sin(\frac{\left|\lambda\right|}{\vartheta}z)+B\cos(\frac{\left|\lambda\right|}{\vartheta}z)$.
To fulfil the boundary condition (\ref{eq:eigen_problem_BC_left}),
$B=0$ must be valid, and hence, the ansatz simplifies to $\phi(z)=A\sin(\frac{\left|\lambda\right|}{\vartheta}z)$.
Furthermore, by means of the boundary condition (\ref{eq:eigen_problem_BC_right}),
we find $\partial_{z}\phi(L)=\frac{\left|\lambda\right|}{\vartheta}A\cos(\frac{\left|\lambda\right|}{\vartheta}L)=0$,
which exhibits infinitely many non-trivial solutions for
\begin{equation}
\left|\lambda_{k}\right|=(k-\frac{1}{2})\frac{\pi}{L}\vartheta\,,\label{eq:eigenvalues}
\end{equation}
with $k=1,2,\ldots$. With regard to the investigation of the set
$\mathcal{S}$, the velocity of the vibrating string is of particular
interest. Consequently, since we deduced the (imaginary) eigenvalues
$\lambda_{k}=\pm i\omega_{k}\vartheta$ with $\omega_{k}=(k-\frac{1}{2})\frac{\pi}{L}$,
the ansatz for the general solution of the velocity can be given according
to\begin{subequations}\label{eq:solution_a_b}
\begin{align}
\tilde{w}_{t}(z,t) & \!=\!\overset{\infty}{\underset{\mathop{k=1}}{\mathop{\sum}}}\!(a_{k}\cos(\omega_{k}\vartheta t)\!+\!b_{k}\sin(\omega_{k}\vartheta t))\varphi_{k}(z),\label{eq:solution_w_t}
\end{align}
where the coefficients $A_{k}$ are hidden in $a_{k}$ and $b_{k}$,
and therefore, for the eigenfunctions we use $\varphi_{k}(z)=\sin(\omega_{k}z)$
here and in the sequel. Hence, an integration of (\ref{eq:solution_w_t})
yields
\begin{equation}
\tilde{w}(z,t)\!=\!\overset{\infty}{\underset{\mathop{k=1}}{\mathop{\sum}}}\!(a_{k}\sin(\omega_{k}\vartheta t)\!-\!b_{k}\cos(\omega_{k}\vartheta t))\frac{\varphi_{k}(z)}{\omega_{k}\vartheta}.\label{eq:solution_w}
\end{equation}
\end{subequations}By means of $\sin(x)=\frac{1}{2i}(e^{ix}-e^{-ix})$
and $\cos(x)=\frac{1}{2}(e^{ix}+e^{-ix})$, after a straightforward
computation we can beneficially rewrite (\ref{eq:solution_a_b}) according
to
\begin{multline*}
\left[\begin{array}{c}
\tilde{w}(z,t)\\
\tilde{w}_{t}(z,t)
\end{array}\right]=\overset{\infty}{\underset{\mathop{k=1}}{\mathop{\sum}}}c_{k}e^{i\omega_{k}\vartheta t}\left[\begin{array}{c}
-i\frac{\varphi_{k}}{\omega_{k}\vartheta}\\
\phi_{k}
\end{array}\right]+\ldots\\
\overset{\infty}{\underset{\mathop{k=1}}{\mathop{\sum}}}c_{-k}e^{-i\omega_{k}\vartheta t}\left[\begin{array}{c}
i\frac{\varphi_{k}}{\omega_{k}\vartheta}\\
\phi_{k}
\end{array}\right]\,,
\end{multline*}
where the coefficients $c_{k}=\frac{1}{2}\left(a_{k}-ib_{k}\right)$
and $c_{-k}=\frac{1}{2}\left(a_{k}+ib_{k}\right)$ fulfil (see \cite[Eq. (3.19)]{Guo2011})
\begin{equation}
\overset{\infty}{\underset{\mathop{k=1}}{\mathop{\sum}}}\left|c_{\pm k}\right|^{2}=\overset{\infty}{\underset{\mathop{k=1}}{\mathop{\sum}}}\left(a_{k}^{2}+b_{k}^{2}\right)<\infty\,,\label{eq:coefficients_bounded}
\end{equation}
which will play an important role later. Thus, we are able to write
$\int_{0}^{L}g(z)\frac{1}{\rho}\tilde{p}\mathrm{d}z=\int_{L_{p_{1}}}^{L_{p_{2}}}\tilde{w}_{t}\mathrm{d}z=0$
as
\begin{equation}
\int_{L_{p_{1}}}^{L_{p_{2}}}\overset{\infty}{\underset{\mathop{k=1}}{\mathop{\sum}}}(c_{k}e^{i\omega_{k}\vartheta t}+c_{-k}e^{-i\omega_{k}\vartheta t})\sin(\omega_{k}z)\mathrm{d}z=0\,.\label{eq:condition_S}
\end{equation}

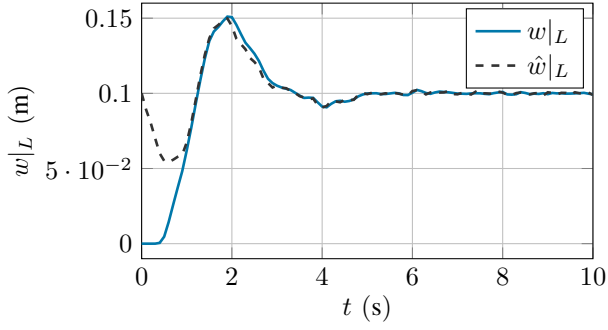
\begin{figure}
\begin{tikzpicture}

\begin{axis}[%
width=6cm,
height=3.4cm,
at={(0.0in,0.0in)},
xlabel style={yshift=0.1cm},
ylabel style={yshift=-0.1cm},
scale only axis,
xmin=0,
xmax=10,
xlabel={$t~(\text{s})$},
xmajorgrids,
ymin=-0.01,
ymax=0.16,
ylabel={$w|_{L}~(\text{m})$},
ymajorgrids,
axis background/.style={fill=white},
legend style={legend cell align=left,align=left,draw=black}
]

\addplot [color=jkuBlue,solid,line width=1.0pt]
  table[row sep=crcr]{%
0	0\\
0.1	1.32096346172423e-13\\
0.2	2.04642351843963e-08\\
0.3	1.14896619622216e-05\\
0.4	0.000512059498051376\\
0.5	0.00466789319270326\\
0.6	0.0144931866424321\\
0.7	0.0260122810061951\\
0.8	0.0372982221074747\\
0.9	0.0482042028291889\\
1	0.0620546887585488\\
1.1	0.0771068380129743\\
1.2	0.0931963247567162\\
1.3	0.108626324483215\\
1.4	0.123069458241493\\
1.5	0.134789996591847\\
1.6	0.141109105259123\\
1.7	0.144130143774684\\
1.8	0.147885690169963\\
1.9	0.151033404739036\\
2	0.150683687338495\\
2.1	0.145566842560852\\
2.2	0.139463699083744\\
2.3	0.133687278609978\\
2.4	0.130149958238749\\
2.5	0.126357462695275\\
2.6	0.121483212573534\\
2.7	0.115081402209043\\
2.8	0.109704011675966\\
2.9	0.106843350887945\\
3	0.105081075002252\\
3.1	0.103907635483112\\
3.2	0.10306471638099\\
3.3	0.10225043359956\\
3.4	0.100398102290234\\
3.5	0.0985608142988808\\
3.6	0.0969605201669917\\
3.7	0.0971557325922409\\
3.8	0.0970059856134233\\
3.9	0.0940276317283797\\
4	0.0912059016372829\\
4.1	0.0907183147898039\\
4.2	0.0930168450483697\\
4.3	0.0947059239567726\\
4.4	0.0950591215229145\\
4.5	0.094419550352293\\
4.6	0.0948990654965872\\
4.7	0.0959400105391529\\
4.8	0.0979824865111642\\
4.9	0.0992098623681997\\
5	0.0999010334054783\\
5.1	0.100007758213401\\
5.2	0.100059480901164\\
5.3	0.0995427586531467\\
5.4	0.0995239026501461\\
5.5	0.0999054025809462\\
5.6	0.100299588439014\\
5.7	0.10000466816442\\
5.8	0.0988971147015873\\
5.9	0.0988431357195226\\
6	0.100424701774005\\
6.1	0.10227647498451\\
6.2	0.101465418713527\\
6.3	0.100233529052875\\
6.4	0.0996752597695701\\
6.5	0.100597186857267\\
6.6	0.101180200332078\\
6.7	0.100906960934783\\
6.8	0.100688314442511\\
6.9	0.100278313118004\\
7	0.0998333472388037\\
7.1	0.0997522201146361\\
7.2	0.0996469638689949\\
7.3	0.0999215439691197\\
7.4	0.100045223985042\\
7.5	0.100018259280398\\
7.6	0.0996578440408692\\
7.7	0.0994262216990313\\
7.8	0.100128183413336\\
7.9	0.100588038155544\\
8	0.10093238273149\\
8.1	0.0992813420658049\\
8.2	0.0991147731461406\\
8.3	0.0994138082224999\\
8.4	0.100785058214696\\
8.5	0.100259685577856\\
8.6	0.0994697442627611\\
8.7	0.099772811799993\\
8.8	0.0995976662094128\\
8.9	0.100080723955278\\
9	0.0994460865538572\\
9.1	0.100249910578145\\
9.2	0.0999800355435301\\
9.3	0.100050900251757\\
9.4	0.0997877023278819\\
9.5	0.0998277257601999\\
9.6	0.0999745423000142\\
9.7	0.100106086807768\\
9.8	0.100343458341899\\
9.9	0.0997548901017225\\
10	0.0993335175965548\\
};

\addlegendentry{$w|_{L}$}
\addplot [color=jkuDarkGrey,dashed,line width=1.0pt]
  table[row sep=crcr]{%
0	0.1\\
0.1	0.0890923509591369\\
0.2	0.0804765879286487\\
0.3	0.070087151586299\\
0.4	0.0600903425768391\\
0.5	0.0550893846860285\\
0.6	0.0543439251108356\\
0.7	0.0558421357307384\\
0.8	0.0578002639165782\\
0.9	0.0592184502037966\\
1	0.0673429101995522\\
1.1	0.080297590246295\\
1.2	0.0958165479426486\\
1.3	0.111667623896279\\
1.4	0.125604548546399\\
1.5	0.136663336050889\\
1.6	0.143147674686778\\
1.7	0.145491767369608\\
1.8	0.148408458453672\\
1.9	0.149805871750661\\
2	0.144767533074895\\
2.1	0.136333246449129\\
2.2	0.132739946154138\\
2.3	0.129330548722276\\
2.4	0.123565032288519\\
2.5	0.120867271696078\\
2.6	0.117179567047812\\
2.7	0.10977390162181\\
2.8	0.10689417796313\\
2.9	0.103890085635254\\
3	0.105182006704655\\
3.1	0.103667618940175\\
3.2	0.103420035440249\\
3.3	0.10169801634955\\
3.4	0.100273186095196\\
3.5	0.0981369271894567\\
3.6	0.0971733037940809\\
3.7	0.0960180360787987\\
3.8	0.0965984767031129\\
3.9	0.0921901921313578\\
4	0.0905601339289446\\
4.1	0.0915938528646395\\
4.2	0.0947920506096003\\
4.3	0.0946203204718879\\
4.4	0.0940710494318604\\
4.5	0.0958235920809382\\
4.6	0.0955028247443654\\
4.7	0.0956888439941983\\
4.8	0.0989812585152426\\
4.9	0.100255115237326\\
5	0.099450899425663\\
5.1	0.100340275594824\\
5.2	0.0997074373471306\\
5.3	0.0993547261904192\\
5.4	0.100016823257072\\
5.5	0.0995703469605144\\
5.6	0.10028827320079\\
5.7	0.0997252060134931\\
5.8	0.0991589697459549\\
5.9	0.0987358771689181\\
6	0.101532277076927\\
6.1	0.102082292238115\\
6.2	0.101240343269493\\
6.3	0.0992883756575382\\
6.4	0.100589830847003\\
6.5	0.100998469759674\\
6.6	0.100510787008915\\
6.7	0.100638010698079\\
6.8	0.101637346921585\\
6.9	0.0995271835796149\\
7	0.0993638228449733\\
7.1	0.100305189556556\\
7.2	0.0997801715736208\\
7.3	0.0992229134019665\\
7.4	0.101078915538025\\
7.5	0.0990272284953368\\
7.6	0.100002613943671\\
7.7	0.099869472979802\\
7.8	0.0995662485355835\\
7.9	0.101211575634084\\
8	0.100467954138067\\
8.1	0.0991238978402305\\
8.2	0.0989336092299525\\
8.3	0.10015483015031\\
8.4	0.100566435033258\\
8.5	0.100014914482967\\
8.6	0.0993015146504244\\
8.7	0.100005792717878\\
8.8	0.100225183859547\\
8.9	0.0985652185217305\\
9	0.100689649324091\\
9.1	0.100310553234478\\
9.2	0.099237083496797\\
9.3	0.100331969553833\\
9.4	0.100267278770196\\
9.5	0.0989626116992087\\
9.6	0.100704807097463\\
9.7	0.100053041795702\\
9.8	0.0996916661524721\\
9.9	0.100489921295775\\
10	0.0988857063634776\\
};
\addlegendentry{$\hat{w}|_{L}$}
\end{axis}
\end{tikzpicture}%
\caption{\label{fig:Comp_w_state_obs}Comparison of the string deflection $w(L,t)$
and the observer state $\hat{w}(L,t)$.}
\end{figure}

Now, we show that the only solution in $\mathcal{S}$ is the trivial
one, i.e. $c_{\pm k}=0\forall k\geq1$ is valid. Otherwise, if there
exists a $k_{0}$ with $\left|c_{k_{0}}\right|\neq0$, due to (\ref{eq:coefficients_bounded})
we can find a $K>k_{0}$ such that\begin{subequations}\label{eq:bound_inf_coeffs}
\begin{align}
\left|\int_{L_{p_{1}}}^{L_{p_{2}}}\overset{\infty}{\underset{\mathop{k=K}}{\mathop{\sum}}}c_{k}\varphi_{k}\mathrm{d}z\right| & <\left|\frac{c_{k_{0}}}{4}\int_{L_{p_{1}}}^{L_{p_{2}}}\varphi_{k_{0}}\mathrm{d}z\right|\\
\left|\int_{L_{p_{1}}}^{L_{p_{2}}}\overset{\infty}{\underset{\mathop{k=K}}{\mathop{\sum}}}c_{-k}\varphi_{k}\mathrm{d}z\right| & <\left|\frac{c_{k_{0}}}{4}\int_{L_{p_{1}}}^{L_{p_{2}}}\varphi_{k_{0}}\mathrm{d}z\right|
\end{align}
\end{subequations}holds, i.e. the sum of the coefficients from $K$
to $\infty$ multiplied with their corresponding eigenfunctions can
be bounded by $c_{k_{0}}$ and $\varphi_{k_{0}}$. Here, it is assumed
that $\int_{L_{p_{1}}}^{L_{p_{2}}}\sin(\omega_{k_{0}})\mathrm{d}z\neq0$,
i.e. $\omega_{k_{0}}\neq\frac{2\pi}{L_{p_{2}}-L_{p_{1}}}j$ with $j\in\mathbb{N}_{+}$.
However, if we consider the absolute value of the eigenvalues (\ref{eq:eigenvalues}),
we find that this is ensured for a proper choice of the length of
the in-domain actuator according to $L_{p_{2}}-L_{p_{1}}\neq\frac{4L}{2k-1}$.
Consequently, because $\omega_{k}\neq\omega_{l}\forall k\neq l$,
for $t>0$ we can reformulate (\ref{eq:condition_S}) as
\begin{multline}
-c_{k_{0}}\int_{L_{p_{1}}}^{L_{p_{2}}}\varphi_{k_{0}}\mathrm{d}z=\int_{L_{p_{1}}}^{L_{p_{2}}}\{\overset{K}{\underset{\mathop{k=1,k\neq k_{0}}}{\mathop{\sum}}}c_{k}e^{i(\omega_{k}-\omega_{k_{0}})\vartheta t}\varphi_{k}+\\
+\overset{\infty}{\underset{\mathop{k=K+1}}{\mathop{\sum}}}c_{k}e^{i(\omega_{k}-\omega_{k_{0}})\vartheta t}\varphi_{k}+\overset{K}{\underset{\mathop{k=1}}{\mathop{\sum}}}c_{-k}e^{-i(\omega_{k}+\omega_{k_{0}})\vartheta t}\varphi_{k}+\\
+\overset{\infty}{\underset{\mathop{k=K+1}}{\mathop{\sum}}}c_{-k}e^{-i(\omega_{k}+\omega_{k_{0}})\vartheta t}\varphi_{k}\}\mathrm{d}z\,.\label{eq:condition_S_ref}
\end{multline}
Next, the idea is to integrate (\ref{eq:condition_S_ref}) with respect
to the time $t$ and to investigate the absolute value. Hence, we
find that the right-hand side of
\begin{multline}
\left|c_{k_{0}}\int_{L_{p_{1}}}^{L_{p_{2}}}\varphi_{k_{0}}\mathrm{d}z\right|t\leq\\
2\left|\int_{0}^{t}\int_{L_{p_{1}}}^{L_{p_{2}}}\{\overset{K}{\underset{\mathop{k=1,k\neq k_{0}}}{\mathop{\sum}}}c_{k}e^{i(\omega_{k}-\omega_{k_{0}})\vartheta\tau}\varphi_{k}\}\mathrm{d}z\mathrm{d}\tau\right|+\\
+2\left|\int_{0}^{t}\int_{L_{p_{1}}}^{L_{p_{2}}}\{\overset{K}{\underset{\mathop{k=1}}{\mathop{\sum}}}c_{-k}e^{-i(\omega_{k}+\omega_{k_{0}})\vartheta\tau}\varphi_{k}\}\mathrm{d}z\mathrm{d}\tau\right|\,,\label{eq:inequaltiy_coefficients}
\end{multline}
where we used (\ref{eq:bound_inf_coeffs}) to obtain an estimation
for the sums from $k=K+1$ to $k=\infty$, is bounded for all $t\geq0$.
Since for an appropriate choice of the actuator-length it is ensured
that the integral on the left-hand side cannot vanish, the only possibility
that inequality (\ref{eq:inequaltiy_coefficients}) holds for $t\rightarrow\infty$
is that $c_{k_{0}}=0$ is valid. Thus, it is shown that the only possible
solution in $\mathcal{S}$ is the trivial one, which finally proves
the asymptotic stability of the observer error and therefore justifies
the application of the observer developed in \cite{Malzer2020}. Furthermore,
in Figure (\ref{fig:Comp_w_state_obs}), the comparison of the string
deflection $w(L,t)$ and the observer state $\hat{w}(L,t)$ is depicted,
where the tip of the string is moved from $w(L,0)=0$ to $w(L,t_{end})=0.1$
and the observer state is initialised as $\hat{w}(L,0)=0.1$.

\section{Conclusion and Outlook}

In this paper, the asymptotic stability of an observer error of an
in-domain actuated vibrating string, where the observer has been developed
in \cite{Malzer2020}, was investigated. First, we showed that the
linear operator, which describes the observer error as an abstract
Cauchy problem, is the infinitesimal generator of a contraction semigroup.
Second, by means of LaSalle's invariance principle the asymptotic
stability of the observer error was proven. In fact, by choosing the
length of the actuator properly, it was shown that the only possible
solution for $\dot{\tilde{\mathscr{H}}}=0$ is the trivial one, which
implies that the observer error tends to zero. Future-research tasks
might deal with the stability analysis of the closed loop obtained
by the controller design presented in \cite{Malzer2020}, or even
with the stability investigation of the combination of controller
and observer.







\bibliographystyle{ieeetr}
\bibliography{my_bib}

\begin{thebibliography}{10}

\bibitem{Luo1998}
Z.-H. Luo, B.-Z. Guo, and O.~Morgul, {\em {Stability and Stabilization of
  Infinite Dimensional Systems with Applications}}.
\newblock {Springer}, 1998.

\bibitem{Miletic2015}
M.~Mileti\'{c}, D.~St{\"u}rzer, and A.~Arnold, ``{An Euler-Bernoulli beam with
  nonlinear damping and a nonlinear spring at the tip},'' {\em {Discrete \&
  Continuous Dynamical Systems}}, vol.~20, no.~9, pp.~3029--3055, 2015.

\bibitem{Stuerzer2016}
D.~St{\"u}rzer, A.~{Arnold}, and A.~{Kugi}, ``{Closed-loop Stability Analysis
  of a Gantry Crane with Heavy Chain},'' {\em {International Journal of
  Control}}, vol.~91, no.~8, pp.~1931--1943, 2018.

\bibitem{Henikl2012}
J.~Henikl, J.~Schr{\"o}ck, T.~Meurer, and A.~Kugi, ``{Infinit-dimensionaler
  Reglerentwurf für Euler-Bernoulli Balken mit Macro-Fibre Composite
  Aktoren},'' {\em {at - Automatisierungstechnik}}, vol.~60, no.~1, pp.~10--19,
  2012.

\bibitem{Ennsbrunner2005}
H.~Ennsbrunner and K.~Schlacher, ``{On the Geometrical Representation and
  Interconnection of Infinite Dimensional Port Controlled Hamiltonian
  Systems},'' {\em {Proceedings of the 44th IEEE Conference on Decision and
  Control and the European Control Conf.}}, no.~5263--5268, 2005.

\bibitem{Schoeberl2012}
M.~Sch{\"o}berl and A.~Siuka, ``{On the port-Hamiltonian representation of
  systems described by partial differential equations},'' in {\em {Proceedings
  of the 4th IFAC Workshop on Lagrangian and Hamiltonian Methods for Non Linear
  Control}}, vol.~45, pp.~1--6, 2012.

\bibitem{Schoeberl2014a}
M.~Sch{\"o}berl and A.~Siuka, ``{Jet bundle formulation of infinite-dimensional
  port-Hamiltonian systems using differential operators},'' {\em {Automatica}},
  vol.~50, no.~2, pp.~607--613, 2014.

\bibitem{Schaft2002}
A.~J. van~der Schaft and B.~Maschke, ``{Hamiltonian formulations of distributed
  parameter systems with boundary energy flow},'' {\em {Journal of Geometry and
  Physics}}, vol.~42, no.~1-2, pp.~166--194, 2002.

\bibitem{Gorrec2005}
Y.~\mbox{Le Gorrec}, H.~J. Zwart, and B.~Maschke, ``Dirac structures and
  boundary control systems associated with skew-symmetric differential
  operators,'' {\em {SIAM J. Control Optim.}}, vol.~44, no.~5, pp.~1864--1892,
  2005.

\bibitem{Schoeberl2013b}
M.~Sch{\"o}berl and A.~Siuka, ``{Analysis and Comparison of Port-Hamiltonian
  Formulations for Field Theories - demonstrated by means of the Mindlin
  plate},'' in {\em {Proceedings of the European Control Conference (ECC)}},
  pp.~548--553, 2013.

\bibitem{Malzer2020}
T.~Malzer, J.~Toledo, Y.~\mbox{Le Gorrec}, and M.~Sch{\"o}berl, ``{Energy-Based
  In-Domain Control and Observer Design for Infinite-Dimensional
  Port-Hamiltonian Systems},'' {\em {arXiv preprint arXiv:2002.01717}}, 2020.

\bibitem{Schoeberl2011}
M.~Sch{\"o}berl and A.~Siuka, ``{On Casimir Functionals for Field Theories in
  Port-Hamiltonian Description for Control Purposes},'' in {\em {Proceedings of
  the 50th IEEE Conference on Decision and Control and European Control
  Conference (CDC-ECC)}}, pp.~7759--7764, 2011.

\bibitem{Rams2017a}
H.~Rams and M.~Sch{\"o}berl, ``{On Structural Invariants in the Energy Based
  Control of Port-Hamiltonian Systems with Second-Order Hamiltonian},'' in {\em
  {Proceedings of the American Control Conference (ACC)}}, pp.~1139--1144,
  2017.

\bibitem{Macchelli2004}
A.~Macchelli and C.~Melchiorri, ``{Modeling and control of the Timoshenko beam.
  The distributed port Hamiltonian approach},'' {\em {SIAM J. Control Optim.}},
  vol.~43, no.~2, pp.~743--767, 2004.

\bibitem{Macchelli2017}
A.~Macchelli, Y.~\mbox{Le Gorrec}, H.~Ramirez, and H.~Zwart, ``{On the
  Synthesis of Boundary Control Laws for Distributed Port-Hamiltonian
  Systems},'' {\em {IEEE Trans. Autom. Control}}, vol.~62, no.~4,
  pp.~1700--1713, 2017.

\bibitem{Toledo2019}
J.~Toledo, H.~Ramirez, Y.~Wu, and Y.~\mbox{Le Gorrec}, ``{Passive observers for
  distributed port-Hamiltonian systems},'' in {\em {InProceedings of the 21st
  IFAC World Congress, Berlin, Germany, July 12-17, 2020 (accepted)}}, 2019.

\bibitem{Schoeberl2008a}
M.~Sch{\"o}berl, H.~Ennsbrunner, and K.~Schlacher, ``{Modelling of
  piezoelectric structures - a Hamiltonian approach},'' {\em {Mathematical and
  Computer Modelling of Dynamical Systems}}, vol.~14, no.~3, pp.~179--193,
  2008.

\bibitem{Schoeberl2014}
M.~Sch{\"o}berl, {\em {Contributions to the Analysis of Structural Properties
  of Dynamical Systems in Control and Systems Theory: A Geometric Approach}}.
\newblock {Shaker Verlag Aachen}, 2014.

\bibitem{Rams2017b}
H.~Rams, M.~Sch{\"o}berl, and K.~Schlacher, ``{Optimal Motion Planning and
  Energy-Based Control of a Single Mast Stacker Crane},'' {\em {IEEE
  Transactions on Control Systems Technology}}, vol.~26, no.~4, pp.~1449--1457,
  2018.

\bibitem{Liu1999}
Z.~Liu and S.~Zheng, {\em {Semigroups Associated with Dissipative Systems}}.
\newblock {Research Notes in Mathematics Series}, Chapman and Hall/CRC, 1999.

\bibitem{Adams2003}
R.~A. Adams and J.~F. Fournier, {\em {Sobolev Spaces}}, vol.~140 of {\em {Pure
  and Applied Mathematics}}.
\newblock {Academic Press, Inc.}, 2nd~ed., 2003.

\bibitem{Guo2011}
W.~Guo and B.-Z. Guo, ``{Parameter estimation and stabilisation for a
  one-dimensional wave equation with boundary output constant disturbance and
  non-collocated control},'' {\em {International Journal of Control}}, vol.~84,
  no.~2, pp.~381--395, 2011.

\end{thebibliography}

\end{document}